\documentclass[12pt,reqno,a4paper]{article}
\usepackage{mathtools}
\usepackage{enumerate}
\usepackage{graphicx} 
\usepackage{amsmath,amssymb,amsfonts,amsthm}
\usepackage{xcolor}
\usepackage[
    inner=2.3cm,
    outer=2.3cm,
    top=2.7cm,
    marginparwidth=2.5cm,
    marginparsep=0.1cm
    ]
    {geometry}

\newtheoremstyle{special}%
{}%
{}%
{}%
{}%
{\scshape}%
{.}%
{.5em}%
{}
\newtheorem{theorem}{Theorem}

\newtheorem{lemma}[theorem]{Lemma}
\newtheorem{proposition}[theorem]{Proposition}
\newtheorem{definition}[theorem]{Definition}
\theoremstyle{definition}
\newtheorem{remark}[theorem]{Remark}

\def\Id{\text{\rm Id}}
\def\Ker{\text{\rm Ker}}

\title{Characterizing nonuniform hyperbolicity by Mather-type admissibility}
\author{Robin Chemnitz\footnote{Mathematics Institute, Freie Universität Berlin, Germany. {\tt E-mail: r.chemnitz@fu-berlin.de}.}\and Davor Dragi\v cevi\'c\footnote{Faculty of Mathematics, University of Rijeka, Rijeka Croatia. {\tt E-mail: ddragicevic@math.uniri.hr}.}}

\begin{document}

\maketitle

\begin{abstract}
    We consider linear cocycles acting on Banach spaces which satisfy the assumptions of the multiplicative ergodic theorem. A cocycle is nonuniformly hyperbolic if all Lyapunov exponents are non-zero, which is equivalent to the existence of a tempered exponential dichotomy. We provide an equivalent characterization of nonuniform hyperbolicity in terms of a Mather-type admissibility of a pair of weighted function spaces. As an application we give a short proof of the robustness of  tempered exponential dichotomies under small linear perturbation.
\end{abstract}

\section{Introduction}
A key feature in the study of (nonautonomous) linear dynamical systems is the existence of (exponentially) stable and (exponentially) unstable subspaces. These subspaces are critical in the description of the long-term evolution of a linear system, governing both qualitative and quantitative behavior. If the state space can be split completely into a stable and an unstable subspace, the dynamic is called \emph{hyperbolic}. Understanding and characterizing hyperbolicity is one of the central topics in the field of dynamical systems. 

In this article, we consider an operator theoretic approach to the (nonuniform) hyperbolicity of a linear cocycle over an ergodic base. A common technique to turn a nonautonomous system into an autonomous one is to augment the state space to include time as an additional dimension. In the case of a cocycle over a base dynamics, this augmentation is particularly useful since we can augment with respect to the base instead of with respect to time. The idea to characterize hyperbolicity using the augmentation of a cocycle over a base dynamics goes back to the work of J.~Mather. We give a brief exposition to these initial results.

Let $f\colon M\to M$ be a $C^1$-diffeomorphism of a compact Riemannian manifold $M$. By $\Gamma (TM)$ we denote the Banach space of all continuous sections $v\colon M\to TM$ equipped with the supremum (maximum) norm $\|v\|_\infty:=\sup_{x\in M}\|v(x)\|$. We define a bounded linear operator (Mather operator) $\mathbb A\colon \Gamma(TM)\to \Gamma (TM)$ by 
\[
[\mathbb A v](x)=Df(f^{-1}(x))v(f^{-1}(x)), \quad x\in M, \ v\in \Gamma (TM).
\]
In his seminal paper~\cite{mather}, J. Mather proved that the following three statements are equivalent:
\begin{enumerate}
\item[$(i)$] $f$ is an Anosov diffeomorphism;
\item[$(ii)$] the spectrum of $\mathbb A$ does not intersect the unit circle $S^1$;
\item[$(iii)$] $1$ is not in the spectrum of $\mathbb A$.
\end{enumerate}
We note that the property $(iii)$ is equivalent to require that for each $w\in \Gamma(TM)$ there exists a unique $v\in \Gamma(TM)$ such that 
\[
v(x)-Df(f^{-1}(x))v(f^{-1}(x))=w(x), \quad x\in M.
\]
This alternative formulation is analogous to the concept of \emph{admissibility} introduced in the context of nonautonomous differential equations by Massera and Sch{\"a}ffer~\cite{massera1966linear} (building on the earlier work of Perron~\cite{Perron} and Li~\cite{Li}), where $w$ is regarded as \emph{input} while $v$ is regarded as an \emph{output}. We refer to~\cite{barreira2018admissibility} for  a survey on admissibility techniques in the study of asymptotic behaviour of nonautonomous dynamics.

Since the original work of Mather there has a substantial progress on extending his results to the case of (not necessarily derivative) linear cocycles acting on arbitrary Banach spaces. More precisely, there are numerous results which characterize \emph{uniform} hyperbolicity of linear cocycles in terms of spectral properties of the Mather operator (in the case of cocycles over maps) or of a generator of the Mather semigroup (in the case of cocycles over flows) acting on suitable Banach spaces. We refer to~\cite{latushkin1999evolution, latushkin1991weighted,  latushkin1995evolutionary, latushkin1996spectral, latushkin1999spectral} and also to~\cite[Chapters 6 and 7]{chicone1999evolution} for a detailed exposition of this theory.

We note that there exists a different but somewhat related approach to the characterization of hyperbolicity for linear cocycles that has been initiated by Chow and Leiva~\cite{chow1995existence}. While in the above mentioned approach by Mather the input and output spaces consist of maps from the base space of our dynamics to the Banach space $X$ on which our cocycle acts, in the second approach  these consist of (two-sided) sequences of elements in $X$ corresponding to a fixed orbit in the base. Hence, an admissibility condition of Mather-type considers all orbits of the base simultaneously, while Chow-Leiva-type admissibility considers each orbit separately. We note that in some cases it is possible to use one approach to tackle the other (see~\cite{latushkin1999spectral}).

Starting from the pioneering works of Oseledets~\cite{Ose} and particularly Pesin~\cite{Pesin} there has been a growing interest in the theory of \emph{nonuniformly} hyperbolic dynamical systems.
The notion of nonuniform hyperbolicity is intimately connected with the non-vanishing of Lyapunov exponents (see Proposition~\ref{criteria}) and represent a far reaching generalization of the notion of uniform hyperbolicity. 
We refer to~\cite{barreira2007nonuniform} for a detailed exposition of this theory. 

Thus, it is natural to ask whether there exists a characterization of nonuniform hyperbolicity in terms of the Mather operator. In this direction, Barreira, Popescu and Valls~\cite{BPV}, as well as Barreira and Valls~\cite{BV} provided characterizations of nonuniform hyperbolicity via hyperbolicity of the Mather operator or of the Mather semigroup. The main feature of these results is that the Mather operator and semigroup act on Banach spaces which are constructed using the so-called Lyapunov norms which transform nonuniform behaviour into the uniform. On the other hand, such norms are difficult to construct before knowing that the dynamics is nonuniformly hyperbolic. Moreover, 
in~\cite{BV} the authors do not consider  the general case of nonuniform hyperbolicity (see Definition~\ref{TED}) but rather require that the random variable ($K$ in Definition~\ref{TED}) measuring the speed of contraction/expansion along stable/unstable direction belongs to some $L^p$ space. 

The main objective of our work is to obtain a Mather type characterization of nonuniform hyperbolicity for linear cocycles over maps without the use of Lyapunov norms. This has a consequence that the input and output spaces do not coincide and that the Mather operator itself does not act as a bounded operator from the input space to the output space. Consequently, in sharp contrast to the standard approach outlined in~\cite{chicone1999evolution} (which is also used in~\cite{BPV}  and~\cite{BV}) we are unable to rely on spectral theory. Instead our arguments are based on the (multiplicative) ergodic theory.
We emphasize that our results hold for strongly measurable cocycles, while the previous work~\cite{BPV} considers only strongly continuous cocycles. Strong measurability is a natural assumption since it is a standard requirement in the infinite-dimensional versions of the Oseledets' theorem  (see Theorem~\ref{MET}) which gives the existence of Lyapunov exponents.

The main results of our work can be summarized in form of the following theorem. For the formal definitions of tempered exponential dichotomies and admissibility, we refer to Section \ref{sec:prelim}.
\begin{theorem}\label{mainThm}
    Let $\mathcal A\colon \Omega \times \mathbb N\to \mathcal L(X)$ be a compact linear cocycle with the property that its generator $A\colon \Omega \to \mathcal L(X)$ satisfies $\log^+ \|A(\cdot)\|\in L^1(\Omega, \mathbb R)$. The following statements are equivalent:
    \begin{enumerate}[(i)]
        \item $\mathcal{A}$ admits a tempered exponential dichotomy;
        \item there exists a random variable $C\colon \Omega \to (0, \infty)$ such that the pair $(L^\infty(\Omega, X), L_C^\infty(\Omega, X))$ is admissible;
        \item there exists a random variable $C\colon \Omega \to (0, \infty)$ such that the pair $(L^\infty_C(\Omega, X), L^\infty(\Omega, X))$ is admissible.
    \end{enumerate}
\end{theorem}
The majority of this article is dedicated to the proof of this theorem. The implications $(i)\Rightarrow (ii)$ and $(i) \Rightarrow (iii)$ are covered by Proposition~\ref{thm:L_LC} and Proposition~\ref{thm:LC_L}, while the implications $(ii) \Rightarrow (i)$ and $(iii) \Rightarrow (i)$ are covered by Theorem~\ref{converse} and Theorem~\ref{converse2}. We note that Theorem~\ref{mainThm} holds more generally for quasi-compact cocycles with an index of compactness of less than $0$  (see Remark \ref{remark}).

We highlight a lemma that is needed in the proofs of Theorem~\ref{converse} and Theorem~\ref{converse2}. To the best of our knowledge a result of this type is new and may be useful in other contexts, see Remark \ref{rem:shear} for a brief discussion of the result.
\theoremstyle{plain}
\newtheorem*{lem_LE}{Lemma \ref{lem:zero_LE_regularity}}
\begin{lem_LE}
    Assume that $(\Omega, \mathcal{F})$ is a Polish space and that $\mathcal A\colon \Omega \times \mathbb N\to \mathcal L(X)$ is a compact linear cocycle satisfying $\log^+\|A(\cdot)\|\in L^1(\Omega, \mathbb R)$ and with  a zero Lyapunov exponent. Then, there is a measurable map $v:\Omega \to X$ such that for $\mathbb P$-a.e.~$\omega \in \Omega$ we have $\|v(\omega)\| = 1$ and 
        \begin{equation*}
            \liminf_{n\to \infty} \|\mathcal{A}(\omega, n) v(\omega)\| \leq 1,  \qquad \limsup_{n\to \infty} \|\mathcal{A}(\omega, n) v(\omega)\| \geq 1.
        \end{equation*}    
\end{lem_LE}

The rest of the article is structured as follows. Section \ref{sec:prelim} describes the setting we work in and contains a version of the multiplicative ergodic theorem as well as the definitions of a tempered exponential dichotomy and admissibility. In Section \ref{sec:dicho_to_adm} we state and prove Proposition~\ref{thm:L_LC} and Proposition~\ref{thm:LC_L} which show that the existence of a tempered dichotomy implies the admissibility of a pair of weighted function spaces. Section \ref{sec:adm_to_dicho} contains the converse results, namely Theorem~\ref{converse} and Theorem~\ref{converse2}. Lastly, in Section \ref{sec:robustness}, we apply our results to show that the existence of a tempered exponential dichotomy is robust under small linear perturbations of the cocycle.

\section{Preliminaries}\label{sec:prelim}
Let $X=(X, \|\cdot \|)$ be an arbitrary separable Banach space. By $\mathcal L(X)$ we will denote the space of all bounded linear operators on $X$ equipped with the operator norm that we will also denote by $\|\cdot \|$. 

Let $(\Omega, \mathcal F, \mathbb P)$ be a probability space such that $(\Omega, \mathcal F)$ is a Polish space equipped with a $\mathbb P$-preserving transformation $\sigma \colon \Omega \to \Omega$ which is invertible with measurable inverse. We will assume that $\mathbb P$ is ergodic with respect to $\sigma$ and that $\sigma$ is aperiodic, i.e.~the set of periodic points is a nullset with respect to $\mathbb P$.

Throughout this paper $\mathbb N$ will denote the set of all nonnegative integers. 
A map $\mathcal A\colon \Omega \times \mathbb N\to \mathcal L(X)$ is said to be a (strongly measurable) \emph{linear cocycle} if the following holds:
\begin{itemize}
\item $\mathcal A(\omega, 0)=\Id$ for $\omega \in \Omega $, where $\Id$ denotes the identity operator on $X$;
\item for $\omega \in \Omega$ and $n, m\in \mathbb N$,
\[
\mathcal A(\omega, n+m)=\mathcal A( \sigma^m \omega, n)\mathcal A(\omega, m);
\]
\item for each $x\in X$, the map $\omega \mapsto A(\omega)x$ is measurable, where $A(\omega):=\mathcal A(\omega, 1)$ is the \emph{generator} of the cocycle $\mathcal A$.
\end{itemize}
We say that $\mathcal{A}$ is a \emph{compact cocycle} if $A(\omega)$ is a compact operator for all $\omega \in \Omega$. The following is a version of the celebrated multiplicative ergodic theorem for compact cocycles.
\begin{theorem}\label{MET}
    Let $\mathcal{A}\colon \Omega \times \mathbb N\to \mathcal L(X)$ be a compact linear cocycle over a Polish space $(\Omega, \mathcal{F})$ such that $\log^+\|A(\cdot)\| \in L^1(\Omega, \mathbb R)$. Then, there is a possibly infinite sequence of numbers $\lambda_1 > \lambda_2 >\hdots >-\infty$ called \emph{Lyapunov exponents}. There is a $\sigma$-invariant subset $\Omega'\subset \Omega$ of full $\mathbb P$-measure such that for each $\lambda_i$ and $\omega \in \Omega'$ there is finite-dimensional  \emph{Oseledets space} $E_i(\omega)$ and an finite-codimensional space $V_i(\omega)$ that together satisfy the following properties for all $\omega \in \Omega'$:
    \begin{itemize}
        \item $X = \left(\oplus_{k=1}^i E_k(\omega) \right) \oplus V_i(\omega)$; 
        \item $A(\omega) E_i(\omega) = E_i(\sigma \omega)$, and $A(\omega) V_i(\omega) \subset V_i(\sigma \omega)$;
        \item $\lim\limits_{n\to \infty} \frac{1}{n} \log(\|\mathcal{A}(\omega, n) v\|)= \lambda_i$, for all non-zero $v \in E_i(\omega)$;
        \item  $\lim\limits_{n\to \infty} \frac{1}{n} \log(\|\mathcal{A}(\omega, n)|_{V_i(\omega)}\|) < \lambda_i$.
    \end{itemize}
\end{theorem}
For a proof, we refer to \cite[Theorem A]{gonzalez2014semi}. In their setting, $(\Omega, \mathcal{F}, \mathbb P)$ is assumed to be a Lebesgue space. Since any regular probability measure on a Polish space results in a Lebesgue space, cf.~\cite[Theorem 2.4.1]{ito1984introduction}, our assumptions suffice. We note that their results hold more generally for quasi-compact cocycles. To avoid additional notation and case distinction, we restrict ourselves to compact cocycles. The results of this article can be extended to quasi-compact cocycles with an index of compactness of less than $0$ without significant changes to the proofs (see Remark \ref{remark}). 
\begin{remark}
The  multiplicative ergodic theorem for matrix cocycles was established by Oseledets~\cite{Ose}. The first versions of Oseledets theorem in the infinite-dimensional setting are due to Ruelle~\cite{ruelle1982characteristic} and Man{\'e}~\cite{mane2006lyapounov} for Hilbert and Banach spaces respectively. For a detailed overview of subsequent  developments we refer to the discussion in ~\cite{gonzalez2014semi}. 
 Finally, we emphasize that many aspects of the theory of nonuniformly hyperbolic dynamical systems were also extended to the infinite-dimensional setting (see~\cite{ blumenthal2017entropy} and references therein).
\end{remark}

We also recall the notion of a tempered exponential dichotomy. Namely,  a random variable $K\colon \Omega \to (0, \infty)$ is said to be \emph{tempered} if 
\[
\lim_{n\to \pm \infty}\frac 1 n \log K(\sigma^n \omega)=0, \quad \text{for $\mathbb P$-a.e. $\omega \in \Omega$.}
\]
The following result is well-known, e.g.~\cite[Proposition 4.3.3]{Arnold1998}.
\begin{proposition}\label{prop}
Let $K\colon \Omega \to (0, \infty)$ be a tempered random variable. Then, for each $\varepsilon>0$ there exists a (tempered) random variable $K_\varepsilon \colon \Omega \to (0, \infty)$ such that 
\[
K(\omega)\le K_\varepsilon (\omega) \quad \text{and} \quad K_\varepsilon (\sigma^n \omega)\le K_\varepsilon (\omega)e^{\varepsilon |n|},
\]
for $\mathbb P$-a.e. $\omega \in \Omega$ and $n\in \mathbb Z$.
\end{proposition}

\begin{definition}\label{TED}
A linear cocycle $\mathcal A\colon \Omega \times \mathbb N \to \mathcal L(X)$ is said to admit a \emph{tempered exponential dichotomy} if there exist a $\sigma$-invariant full measure set $\Omega'\subset \Omega$, a tempered random variable $K\colon \Omega \to (0, \infty)$, $\lambda >0$ and a family of projections $\Pi^s(\omega)$ on $X$ such that the following holds for all $\omega \in \Omega'$ and $n\in \mathbb N$:
\begin{itemize}
\item $\Pi^s(\sigma \omega)A(\omega)=A(\omega)\Pi^s(\omega)$;
\item $A(\omega)\rvert_{\Ker \Pi^s(\omega)}\colon \Ker \Pi^s(\omega)\to \Ker \Pi^s(\sigma \omega)$ is invertible and
\begin{align}\label{dic1}
    \|\mathcal A(\omega, n)\Pi^s(\omega)\| &\le K(\omega)e^{-\lambda n}\\\label{dic2}
    \|\mathcal A(\omega, -n)(\Id-\Pi^s(\omega))\|&\le K(\omega)e^{-\lambda n},
\end{align}
where 
\[
\mathcal A(\omega, -n):=\left (\mathcal A(\sigma^{-n} \omega, n)\rvert_{\Ker \Pi^s(\sigma^{-n} \omega)}\right )^{-1}\colon \Ker \Pi^s (\omega)\to \Ker \Pi^s (\sigma^{-n} \omega);
\]
\item for each $v\in X$, $\omega \mapsto \Pi^s(\omega)v$ is measurable.
\end{itemize}
\end{definition}
A tempered exponential dichotomy with $K(\omega) = K>0$ constant and $\Omega' = \Omega$ is called a \emph{uniform exponential dichotomy}, which is a significantly stronger notion. The existence of a tempered exponential dichotomy is closely related to the Lyapunov exponents of the cocycle. The following result follows from~\cite[Proposition 3.2]{backes2019periodic}.
\begin{proposition}\label{criteria}
    Let $\mathcal{A}\colon \Omega \times \mathbb N\to \mathcal L(X)$ be a compact cocycle over a Polish space $(\Omega, \mathcal{F})$ such that $\log^+\|A(\cdot)\| \in L^1(\Omega, \mathbb R)$. The following assertions are equivalent:
    \begin{enumerate}
        \item[(i)] $\mathcal A$ admits a tempered exponential dichotomy;
        \item[(ii)] all Lyapunov exponents of $\mathcal A$ are nonzero.
    \end{enumerate}
\end{proposition}

We also introduce the notion of admissibility. 
\begin{definition}
Let $\mathcal A\colon \Omega \times \mathbb N\to \mathcal L(X)$ be a linear cocycle and let $\mathcal Y_i$, $i=1, 2$ be two Banach spaces consisting of measurable functions from $\Omega$ to $X$. We say that the pair $(\mathcal Y_1, \mathcal Y_2)$ is \emph{admissible} for $\mathcal A$ if for each $g\in \mathcal Y_2$ there exists a unique $f\in \mathcal Y_1$ such that
\begin{equation}\label{adm}
    f(\omega)-A(\sigma^{-1}\omega)f(\sigma^{-1}\omega)=g(\omega), \quad \text{for $\mathbb P$-a.e. $\omega \in \Omega$.}
\end{equation}
\end{definition}

The admissibility condition \eqref{adm} can be formulated in terms of the \emph{Mather-operator}, also called \emph{Mather-semigroup}. The Mather-operator $\mathbb A$ is defined for all measurable functions $f:\Omega \to X$ by
\begin{equation*}
    [\mathbb A f](\omega) = A(\sigma^{-1} \omega) f(\sigma^{-1}\omega).
\end{equation*}
If $(\mathcal Y_1, \mathcal Y_2)$ is an admissible pair, then for every $g\in \mathcal{Y}_2$ there is an $f \in \mathcal{Y}_1$ such that $(\Id - \mathbb A) f = g$. The map $R:\mathcal{Y}_2 \to \mathcal{Y}_1$ that maps each $g\in \mathcal{Y}_2$ to the respective $f \in \mathcal{Y}_1$ can be interpreted as a type of partial inverse in the sense that
\begin{equation*}
    (\Id - \mathbb A) R = \Id_{\mathcal{Y}_2}.
\end{equation*}
However, in general $R (\Id - \mathbb A)$ is not a well-defined operator on $\mathcal{Y}_1$ since $(\Id - \mathbb A)$ does not map to $\mathcal{Y}_2$. 

An admissible pair $(\mathcal Y_1, \mathcal Y_2)$ is often of the type $\mathcal{Y}_2 \subset \mathcal{Y}_1$, e.g.~when the norm of $\mathcal{Y}_2$ is stronger than the norm $\mathcal{Y}_1$. Even though it is not explicitly stated, all admissible pairs we work with in the rest of this article are of this form. In that case, the inverse operator $R$ makes functions $g\in \mathcal{Y}_2$ more irregular. Hence, $R$ is only defined on a space $\mathcal{Y}_2$ of high enough regularity such that $Rg \in \mathcal{Y}_1$ is well-defined.

In many cases $(\mathcal{Y}, \mathcal{Y})$ admissibility is equivalent to the existence of a uniform exponential dichotomy, see for example \cite{latushkin1999evolution} for the case $\mathcal{Y}= C(\Omega, X)$ and \cite{latushkin1991weighted} for the case $\mathcal{Y}=L^2(\Omega, X)$. For a detailed exposition to the Mather-semigroup, its spectral properties and its relations to uniform exponential dichotomies we refer to \cite[Chapters 6 and 7]{chicone1999evolution}.

\section{Tempered dichotomy implies admissibility}\label{sec:dicho_to_adm}
We introduce a class of function spaces. For a random variable $C\colon \Omega \to (0, \infty)$, by $L_C^\infty(\Omega, X)$ we will denote the space of all measurable maps $f\colon \Omega \to X$ such that 
\[
\|f\|_{\infty, C}:=\text{esssup}_{\omega \in \Omega}\left(C(\omega)\|f(\omega)\|\right )<\infty,
\]
where $\text{esssup}$ is taken with respect to $\mathbb P$.
It is straightforward to verify that $(L_C^\infty(\Omega, X), \|~\cdot~\|_{\infty, C})$ is a Banach space after identifying functions which coincide $\mathbb P$-a.e. In the particular case when $C\equiv 1$, we will write $L^\infty(\Omega, X)$ instead of $L_C^\infty(\Omega, X)$ and $\|\cdot \|_\infty$ instead of $\|\cdot \|_{\infty, C}$.

The following are our first results.

\begin{proposition}\label{thm:L_LC}
Let $\mathcal A\colon \Omega \times \mathbb N\to \mathcal L(X)$ be a linear cocycle that admits a tempered exponential dichotomy. Then, there exists a tempered random variable $C\colon \Omega \to (0, \infty)$ with the property that the pair $(L^\infty(\Omega, X), L_C^\infty(\Omega, X))$ is admissible for $\mathcal A$.
\end{proposition}

\begin{proposition}\label{thm:LC_L}
Let $\mathcal A\colon \Omega \times \mathbb N\to \mathcal L(X)$ be a linear cocycle that admits a tempered exponential dichotomy. Then, there exists a tempered random variable $C\colon \Omega \to (0, \infty)$ with the property that
the pair $(L_C^\infty(\Omega, X), L^\infty(\Omega, X))$ is admissible for $\mathcal A$.

\end{proposition}

\begin{proof}[Proof of Proposition~\ref{thm:L_LC}]
Let $K\colon \Omega \to (0, \infty)$ and $\lambda>0$ be as in Definition~\ref{TED}. We choose $C=K$. Take $g\in L_C^\infty(\Omega, X)$ and set
\begin{equation}\label{x}
\begin{split}
f(\omega) &:=\sum_{n=0}^\infty \mathcal A(\sigma^{-n} \omega, n)\Pi^s(\sigma^{-n}\omega)g(\sigma^{-n}\omega)-\sum_{n=1}^\infty \mathcal A(\sigma^n \omega, -n)\Pi^u(\sigma^n \omega)g(\sigma^n \omega),
\end{split}
\end{equation}
for $\omega \in \Omega$, where $\Pi^u(\omega)=\Id-\Pi^s(\omega)$. By~\eqref{dic1} and~\eqref{dic2} we have that 
\[
\begin{split}
\|f(\omega) \| &\le \sum_{n=0}^\infty \|\mathcal A(\sigma^{-n} \omega, n)\Pi^s(\sigma^{-n}\omega)g(\sigma^{-n}\omega)\|
+\sum_{n=1}^\infty \|\mathcal A(\sigma^n \omega, -n)\Pi^u(\sigma^n \omega)g(\sigma^n \omega) \| \\
&\le \sum_{n=0}^\infty K(\sigma^{-n} \omega)e^{-\lambda n}\|g(\sigma^{-n}\omega)\|+\sum_{n=1}^\infty K(\sigma^n \omega)e^{-\lambda n}\|g(\sigma^n \omega)\| \\
&\le \left (\sum_{n=0}^\infty e^{-\lambda n}+\sum_{n=1}^\infty e^{-\lambda n} \right) \|g\|_{\infty, K} \\
&=\frac{1+e^{-\lambda}}{1-e^{-\lambda}}\|g\|_{\infty, K},
\end{split}
\]
for $\mathbb P$-a.e. $\omega \in \Omega$. Hence, $f\in L^\infty(\Omega, X)$.  On the other hand, for $\mathbb P$-a.e. $\omega \in \Omega$ we have that 
 \begin{equation}\label{computation}
 \begin{split}
     f(\omega)-A(\sigma^{-1}\omega)f(\sigma^{-1}\omega) &=\sum_{n=0}^\infty \mathcal A(\sigma^{-n} \omega, n)\Pi^s(\sigma^{-n}\omega)g(\sigma^{-n}\omega) \\
     &\phantom{=}-\sum_{n=0}^\infty \mathcal A(\sigma^{-(n+1)}\omega, n+1)\Pi^s(\sigma^{-(n+1)}\omega)g(\sigma^{-(n+1)}\omega)\\
     &\phantom{=}-\sum_{n=1}^\infty \mathcal A(\sigma^n \omega, -n)\Pi^u(\sigma^n \omega)g(\sigma^n \omega)\\
     &\phantom{=}+\sum_{n=1}^\infty \mathcal A(\sigma^{n-1}\omega, -(n-1))\Pi^u(\sigma^{n-1}\omega)g(\sigma^{n-1}\omega)\\
     &=\sum_{n=0}^\infty \mathcal A(\sigma^{-n} \omega, n)\Pi^s(\sigma^{-n}\omega)g(\sigma^{-n}\omega) \\
     &\phantom{=}-\sum_{n=1}^\infty \mathcal A(\sigma^{-n}\omega, n)\Pi^s(\sigma^{-n}\omega)g(\sigma^{-n}\omega)\\
      &\phantom{=}-\sum_{n=1}^\infty \mathcal A(\sigma^n \omega, -n)\Pi^u(\sigma^n \omega)g(\sigma^n \omega)\\
     &\phantom{=}+\sum_{n=0}^\infty \mathcal A(\sigma^n\omega, -n)\Pi^u(\sigma^{n}\omega)g(\sigma^{n}\omega)\\
     &=\Pi^s(\omega)g(\omega)+\Pi^u(\omega)g(\omega)\\
     &=g(\omega),
 \end{split}
 \end{equation}
 yielding~\eqref{adm}.
 
It remains to establish the uniqueness of $f$. To this end, it is sufficient to show for $f\in L^\infty(\Omega, X)$ such that 
 \begin{equation}\label{f}
 f(\omega)=A(\sigma^{-1}\omega)f(\sigma^{-1}\omega) \quad \text{for $\mathbb P$-a.e. $\omega \in \Omega$,}
 \end{equation}
we have that $f=0$. Let $f^s(\omega):=\Pi^s(\omega)f(\omega)$. Then, \eqref{dic1} and Proposition~\ref{prop} imply that 
\[
\begin{split}
\|f^s(\omega)\| =\|\mathcal A(\sigma^{-n}\omega, n)f^s(\sigma^{-n}\omega)\| &\le K(\sigma^{-n}\omega)e^{-\lambda n}\|f(\sigma^{-n}\omega)\| \\
&\le K_{\lambda/2}(\sigma^{-n}\omega)e^{-\lambda n}\|f(\sigma^{-n}\omega)\| \\
&\le K_{\lambda/2}(\omega)e^{-\frac{\lambda}{2}n}\|f\|_\infty,
\end{split}
\]
for $\mathbb P$-a.e. $\omega \in \Omega$ and $n\in \mathbb N$. Letting $n\to \infty$  yields  that  $f^s(\omega)=0$ for $\mathbb P$-a.e. $\omega \in \Omega$. Similarly one can show that $f^u(\omega)=\Pi^u(\omega)f(\omega)=0$ for $\mathbb P$-a.e. $\omega \in \Omega$, and thus $f=0$. The proof of the theorem is completed.

\end{proof}

The proof of Proposition~\ref{thm:LC_L} is very similar to that of Proposition~\ref{thm:L_LC}.

\begin{proof}[Proof of Proposition~\ref{thm:LC_L}.]
Let $K:\Omega \to (0,\infty)$ and $\lambda>0$ be as in Definition~\ref{TED}. 
Take $g\in L^\infty(\Omega, X)$ and let $f$ be as in~\eqref{x}.
 By~\eqref{dic1}, \eqref{dic2} and Proposition~\ref{prop} we have that 
\[
\begin{split}
\|f(\omega) \| &\le \sum_{n=0}^\infty \|\mathcal A(\sigma^{-n} \omega, n)\Pi^s(\sigma^{-n}\omega)g(\sigma^{-n}\omega)\|
+\sum_{n=1}^\infty \|\mathcal A(\sigma^n \omega, -n)\Pi^u(\sigma^n \omega)g(\sigma^n \omega) \| \\
&\le \sum_{n=0}^\infty K(\sigma^{-n} \omega)e^{-\lambda n}\|g(\sigma^{-n}\omega)\|+\sum_{n=1}^\infty K(\sigma^n \omega)e^{-\lambda n}\|g(\sigma^n \omega)\| \\
&\le \sum_{n=0}^\infty K_{\lambda/3}(\sigma^{-n} \omega)e^{-\lambda n}\|g(\sigma^{-n}\omega)\|+\sum_{n=1}^\infty K_{\lambda/3}(\sigma^n \omega)e^{-\lambda n}\|g(\sigma^n \omega)\| \\
&\le K_{\lambda/3}(\omega) \|g\|_\infty \left (\sum_{n=0}^\infty e^{-\frac{2\lambda}{3}n}+\sum_{n=1}^\infty e^{-\frac{2\lambda}{3}n}\right )\\
&=\frac{1}{C(\omega)}\frac{1+e^{-\frac{2\lambda}{3}}}{1-e^{-\frac{2\lambda}{3}}}\|g\|_\infty,
\end{split}
\]
for $\mathbb P$-a.e. $\omega \in \Omega$, where $C\colon \Omega \to (0, \infty)$ is given by 
\[
C(\omega):=\frac{1}{K_{\lambda/3}(\omega)}.
\]
 We conclude that $f\in L_C^\infty(\Omega, X)$.  Since $K_{\lambda/3}$ is tempered, we have that $C$ is also tempered. The same argument as in the proof of Proposition~\ref{thm:L_LC} shows that~\eqref{adm} holds. 

 In order to establish the uniqueness of $f$, suppose that $f\in L_C^\infty(\Omega, X)$ satisfies~\eqref{f}.  Using the same notation as in the proof of Proposition~\ref{thm:L_LC} we have that
\[
\begin{split}
\|f^s(\omega)\| =\|\mathcal A(\sigma^{-n}\omega, n)f^s(\sigma^{-n}\omega)\| &\le K(\sigma^{-n}\omega)e^{-\lambda n}\|f(\sigma^{-n}\omega)\| \\
&\le (K_{\lambda/3}(\sigma^{-n}\omega))^2 e^{-\lambda n}\|f\|_{\infty, C} \\
&\le  (K_{\lambda/3}(\omega))^2 e^{-\frac{\lambda}{3}n}\|f\|_{\infty, C},
\end{split}
\]
for $\mathbb P$-a.e. $\omega \in \Omega$ and $n\in \mathbb N$. Letting $n\to \infty$ yields that $f^s(\omega)=0$ for $\mathbb P$-a.e. $\omega \in \Omega$. Similarly, $f^u(\omega)=0$ for $\mathbb P$-a.e. $\omega \in \Omega$ and thus $f=0$.

\end{proof}

\section{From admissibility to tempered dichotomy}\label{sec:adm_to_dicho}
We are now interested in formulating converse results to Proposition~\ref{thm:L_LC} and Proposition~\ref{thm:LC_L}. For this, we will restrict our attention to the case of compact cocycles.

\begin{theorem}\label{converse}
Let $\mathcal A\colon \Omega \times \mathbb N\to \mathcal L(X)$ be a compact linear cocycle with the property that its generator $A\colon \Omega \to \mathcal L(X)$ satisfies $\log^+ \|A(\cdot)\|\in L^1(\Omega, \mathbb R)$. Moreover, suppose that there exists a random variable $C\colon \Omega \to (0, \infty)$ with the property that the pair $(L^\infty(\Omega, X), L_C^\infty(\Omega, X))$ is admissible. Then, $\mathcal A$ admits a tempered exponential dichotomy. 
\end{theorem}

\begin{theorem}\label{converse2}
Let $\mathcal A\colon \Omega \times \mathbb N\to \mathcal L(X)$ be a compact linear cocycle with the property that its generator $A\colon \Omega \to \mathcal L(X)$ satisfies $\log^+ \|A(\cdot)\|\in L^1(\Omega, \mathbb R)$. Moreover, suppose that there exists a random variable $C\colon \Omega \to (0, \infty)$ with the property that the pair $(L^\infty_C(\Omega, X), L^\infty(\Omega, X))$ is admissible. Then, $\mathcal A$ admits a tempered exponential dichotomy. 
\end{theorem}

\begin{remark}
We note that in the statement of Theorem~\ref{converse} and Theorem~\ref{converse2} it is not required that $C$ is a tempered random variable.
\end{remark}

We will first establish several auxiliary results.

\begin{lemma}\label{lem:boundedness}
Let $\mathcal A\colon \Omega \times \mathbb N\to \mathcal L(X)$ be a  linear cocycle. The following hold:
\begin{enumerate}
\item[(i)] if the pair $(L^\infty(\Omega, X), L_C^\infty(\Omega, X))$ is admissible for $\mathcal A$ (for some random variable $C\colon \Omega \to (0, \infty)$), then 
there is a constant $L>0$ such that for each $g\in L_C^\infty(\Omega, X)$ we have 
\begin{equation}\label{eq:L}
    \|f\|_{\infty} \le L\|g\|_{\infty, C},
\end{equation}
where $f\in L^\infty(\Omega,  X)$ is such that~\eqref{adm} holds;
\item[(ii)] if the pair $(L_C^\infty(\Omega, X), L^\infty(\Omega, X))$ is admissible for $\mathcal A$ (for some random variable $C\colon \Omega \to (0, \infty)$), then there is a constant $L>0$ such that for each $g\in L^\infty(\Omega, X)$ we have 
\begin{equation}\label{eq:L1}
    \|f\|_{\infty, C} \le L\|g\|_\infty,
\end{equation}
where $f\in L_C^\infty(\Omega,  X)$ is such that~\eqref{adm} holds.
\end{enumerate}
\end{lemma}

\begin{proof}
We will  only establish $(i)$ as $(ii)$ can be proven using the same argument.
We consider the linear operator $R \colon L_C^\infty(\Omega, X)\to L^\infty(\Omega, X)$ defined by $R g=f$, where $f$ is the unique element of $L^\infty(\Omega, X)$ such that~\eqref{adm} holds. We claim that $R$ is a closed operator. To this end, suppose that $(g_n)_{n\in \mathbb N}\subset L_C^\infty(\Omega, X)$ is a sequence such that $g_n\to g$ in $L_C^\infty(\Omega, X)$ and $f_n:=R y_n \to f$ in $L^\infty(\Omega, X)$. This easily yields that 
\begin{equation}\label{xy}
  \lim_{n\to \infty}f_n(\omega)=f(\omega) \quad \text{and} \quad \lim_{n\to \infty}g_n(\omega)=g(\omega), \quad \text{for $\mathbb P$-a.e. $\omega \in \Omega$.} 
\end{equation}
On the other hand, since $f_n=R g_n$ we have that 
\begin{equation}\label{xy1}
    f_n(\omega)-A(\sigma^{-1}\omega)f_{n-1}(\omega)=g_n(\omega), \quad \text{for $\mathbb P$-a.e. $\omega \in \Omega$ and $n\in \mathbb N$.}
\end{equation}
By passing to the limit when $n\to \infty$ in~\eqref{xy1} and using~\eqref{xy} we obtain that 
\[
f(\omega)-A(\sigma^{-1}\omega)f(\omega)=g(\omega), \quad \text{for $\mathbb P$-a.e. $\omega \in \Omega$.}
\]
Therefore $f=R g$ and we conclude that $R$ is a closed operator. It follows from the Closed Graph Theorem that $R$ is bounded, and consequently~\eqref{eq:L} holds with $L:=\lVert R\rVert$.
\end{proof}

\begin{lemma}\label{lem:zero_Birkhoff}
    Let $(\Omega, \mathcal{F}, \mathbb P, \sigma)$ be an ergodic system and let 
    $\varphi \in L^1(\Omega, \mathbb R)$ with $\int_\Omega \varphi\, d\mathbb P = 0$. Then, for $\mathbb P$-a.e.~$\omega \in \Omega$ the Birkhoff sums $S_n\varphi (\omega):=\sum_{k=0}^{n-1} \varphi( \sigma^k(\omega))$ satisfy 
    \begin{equation*}
        \liminf_{n\to \infty} S_n \varphi(\omega) \leq 0 \quad \text{and} \quad  \limsup_{n\to \infty} S_n \varphi(\omega) \geq 0.
    \end{equation*}
\end{lemma}
\begin{proof}
    It suffices to prove the inequality for the $\liminf$ since the inequality for the $\limsup$ follows by considering $-\varphi$ instead of $\varphi$. Assume for contradiction that there is an $\varepsilon>0$ such that the set
    \begin{equation*}
        \Omega_\varepsilon := \left \{ \omega \in \Omega: \liminf_{n\to \infty} S_n\varphi(\omega) > \varepsilon \right \},
    \end{equation*}
    has positive $\mathbb P$-measure. For each $\omega \in \Omega_\varepsilon$, let $n(\omega)$ be the largest positive integers so that $S_n\varphi(\omega)<\frac{\varepsilon}{2}$.  If there is no such integer we set $n(\omega):=0$.
    Then, the point $\omega' := \sigma^{n(\omega)} \omega$ satisfies $S_n\varphi(\omega')\geq c$
    for all $n\geq 1$ with $c := \frac{\varepsilon}{2} -  S_{n(\omega)}\varphi(\omega)>0$. Therefore, the set
    \begin{equation*}
        \Omega^* := \{\omega \in \Omega: \text{$\exists c>0$ such that   $S_n\varphi(\omega) \geq c$   for  every $n \geq 1$}\},
    \end{equation*}
    has positive $\mathbb P$-measure as well. By continuity of $\mathbb P$, there must be some $c>0$ such that the set 
    \begin{equation*}
        \Omega^*_c := \{\omega \in \Omega :  S_n\varphi(\omega) > c, \ \forall n \geq 1\},
    \end{equation*}
    has positive $\mathbb P$-measure. By the Birkhoff ergodic theorem we find
    \begin{equation*}
        \lim_{n\to \infty}\frac{1}{n} \sum_{k=0}^{n-1} \mathbf{1}_{\Omega^*_c}(\sigma^k \omega) = \mathbb P(\Omega^*_c) >0,
    \end{equation*}
    for $\mathbb P$-a.e.~$\omega \in \Omega$, where $\mathbf{1}_{\Omega^*_c}$ denotes the indicator function of the set $\Omega^*_c$.
    For fixed $\omega \in \Omega$, by the definition of $\Omega^*_c$, for every $n\in \mathbb N$ such that  $\sigma^n \omega \in \Omega_c^*$, we find  that
    \[
    S_{n+m}\varphi(\omega)=S_n \varphi(\omega)+S_m \varphi(\sigma^n \omega)>S_n\varphi(\omega)+c, \quad \forall m \geq 1.
    \]
This implies that 
    \begin{equation}\label{eq:zero_birkhoff_contradiction}
        \lim_{n\to \infty} \frac{1}{n} S_n\varphi(\omega) \ge  \lim_{n\to \infty} \frac{1}{n} c \sum_{k=0}^{n-1} 1_{\Omega^*_c}(\sigma^k \omega) = c \mathbb P (\Omega^*_c) >0,
    \end{equation}
    for $\mathbb P$-a.e. $\omega \in \Omega$.
    This is a contradiction to the assumption $\int_\Omega \varphi\, d\mathbb P= 0$, since the Birkhoff ergodic theorem asserts that the limit on the left-hand side in~\eqref{eq:zero_birkhoff_contradiction} is equal to $0$.
\end{proof}

\begin{lemma}\label{lem:zero_LE_regularity} 
    Assume that $(\Omega, \mathcal{F})$ is a Polish space and that $\mathcal A\colon \Omega \times \mathbb N\to \mathcal L(X)$ is a compact linear cocycle satisfying $\log^+\|A(\cdot)\|\in L^1(\Omega, \mathbb R)$ and with   a zero Lyapunov exponent. Then, there is a measurable map $v:\Omega \to X$ such that for $\mathbb P$-a.e.~$\omega \in \Omega$ we have $\|v(\omega)\| = 1$ and 
    \begin{equation}\label{eq:statement_lemma}
        \liminf_{n\to \infty} \|\mathcal{A}(\omega, n) v(\omega)\| \leq 1,  \qquad \limsup_{n\to \infty} \|\mathcal{A}(\omega, n) v(\omega)\| \geq 1.
    \end{equation}
\end{lemma}

\begin{proof}
    For a detailed introduction to invariant/ergodic measures of the skew-product, we refer to~\cite[Section 1.3.4]{kuksin2012mathematics} or \cite[Section 3]{arnold1999jordan}. 
    Consider the projectivized skew-product $\Xi:\Omega \times P(X) \to \Omega \times P(X)$, defined by $\Xi(\omega, \overline{v}) = (\sigma \omega, \overline{A(\omega) v})$, where $\overline{v} \in  P(X)$ denotes projectivization of $v\in X$. We also define the observable $\varphi(\omega, \overline{v}) = \log (\|A(\omega) v\|)$. The observable $\varphi$ is defined so that 
    \begin{equation}\label{eq:LE_property_f}
        \log(\|\mathcal{A}(\omega, n) v\|) = \sum_{k=0}^{n-1} \varphi(\Xi^k (\omega, \overline{v})), 
    \end{equation}
    for $v\in X$ with $\| v\| = 1$.
    
    Each ergodic measure $\mu$ of $\Xi$ has a marginal on $\Omega$ which is ergodic with respect to $\sigma$. An ergodic measure $\mu$ with marginal $\mathbb P$ can be decomposed into sample measures $\mu_\omega$ on $P(X)$ such that
    \begin{equation*}
        \mu(d\omega, dx) =  \mu_\omega(dx) \mathbb P(d\omega).
    \end{equation*}
    The samples measures are $\mathbb P$-a.s.~unique and satisfy $\mu_{\sigma \omega} = \overline{A(\omega)}_* \mu_\omega$, where $\overline{A(\omega)}_*$ is the push-forward  on $P(X)$, see \cite[Proposition 1.3.27]{kuksin2012mathematics}.
    The Oseledets space $E_0$ corresponding to the zero Lyapunov exponent can be identified with a compact random invariant subset of $P(X)$. By \cite[Theorem 1.6.13]{Arnold1998}, there is at least one ergodic measure $\mu$ of $\Xi$ that is supported on $E_0$ in the sense that $\mu_\omega(E_0(\omega))=1$ for $\mathbb P$-a.e.~$\omega \in \Omega$. By \eqref{eq:LE_property_f}, this ergodic measure has to satisfy $\mu(\varphi) = 0$. 
    By Lemma \ref{lem:zero_Birkhoff}, for $\mu$-a.e.~pair $(\omega, \overline{v}) \in \Omega \times P(X)$ we find that
    \[\liminf_{n\to \infty} \sum_{k=0}^{n-1} \varphi(\Xi^k (\omega, \overline{v})) \leq 0, \qquad \limsup_{n\to \infty} \sum_{k=0}^{n-1} \varphi(\Xi^k (\omega, \overline{v})) \geq 0.\]
    In view of \eqref{eq:LE_property_f}, this implies
    \begin{equation}\label{eq:limsup_liminf}
        \liminf_{n\to \infty} \|\mathcal{A}(\omega, n) v\| \leq 1,  \qquad \limsup_{n\to \infty} \|\mathcal{A}(\omega, n) v\| \geq 1,
    \end{equation}
    where $v\in  X$ is normalized.
    
    It remains to show that $v(\omega)$ can be chosen measurably such that \eqref{eq:limsup_liminf} holds. Consider the set
    \begin{equation*}
        F:=\left \{(\omega, \overline v):   \liminf_{n\to \infty} \|\mathcal{A}(\omega, n) v\| \leq 1 \quad \textup{and} \quad  \limsup_{n\to \infty} \|\mathcal{A}(\omega, n) v\| \geq 1\right \} \subset \Omega \times P(X).
    \end{equation*}
    This set is measurable in $\mathcal{F}\otimes \mathcal{B}(P(X))$, where $\mathcal{B}(P(X))$ denotes the Borel $\sigma$-algebra on $P(X)$.
    Since $\Omega$ and $X$ (and thereby $P(X)$) are assumed to be Polish spaces, $F$ is Borel-measurable and, hence, analytic, cf.~\cite[Proposition 8.2.3]{cohn2013measure}. We know that for $\mathbb P$-a.e.~$\omega \in \Omega$ the section $F_\omega:= F\cap \big(\{\omega\} \times P(X)\big)$ is non-empty. Therefore, the projection $\pi_\Omega(F):= \{\omega: F_\omega \neq \emptyset\}$ is measurable with respect to the $\mathbb P$-completion $\mathcal{F}_*$ and has full $\mathbb P$-measure. Since $F$ is analytic, \cite[Theorem 8.5.3]{cohn2013measure} states that there is a $\mathcal{F}_U-\mathcal{B}(P(X))$-measurable map $v_*:\pi_\Omega(F) \to P(X)$ with $v_*(\omega) \in F_\omega$ for all $\omega \in \pi_\Omega(F)$, where $\mathcal{F}_U$ is the universal completion of $\mathcal{F}$. Since $\mathcal{F}_U \subset \mathcal{F}_*$, we conclude that $v_*:\pi_\Omega(F) \to P(X)$ is $\mathcal{F}_*-\mathcal{B}(P(X))$-measurable. Since $\pi_\Omega(F)$ has full $\mathbb P$-measure, the map $v_*$ can be extended to all of $\Omega$ while remaining $\mathcal{F}_*-\mathcal{B}(P(X))$-measurable. Now, \cite[Lemma 1.2]{crauel2002random} states that there is a $\mathcal{F}-\mathcal{B}(P(X))$-measurable map $v:\Omega \to P(X)$ that coincides with $v_*$ for $\mathbb P$-a.e.~$\omega \in\Omega$. This map $v$ still satisfies $v(\omega) \in F_\omega$ for $\mathbb P$-a.e.~$\omega \in \Omega$. 
    
    It remains to concatenate $v$ with a measurable map $s:P(X) \to X$ such that $\overline{s(\overline{v})} = \overline{v}$, i.e.~$s$ measurably chooses one of the two possible representations of a projective vector $\overline{v}$. The existence of a map follows from e.g.~\cite[Theorem 2.6]{crauel2002random}. This completes the proof.
\end{proof}
\begin{remark}\label{rem:shear}
    We emphasize that in general the vectors in the Oseledets space corresponding to a zero Lyapunov exponent do not satisfy \eqref{eq:statement_lemma}. While it is guaranteed that the norm of any such vector grows at most subexponentially, it may still grow polynomially. As a simple illustrative example,  we consider the cocycle with constant generator
    \begin{equation*}
        A(\omega) = \begin{pmatrix}
            1 & 1 \\ 0 & 1
        \end{pmatrix}.
    \end{equation*}
    It is easily verified that the entirety of $\mathbb R^2$ is the Oseledets space corresponding to a zero Lyapunov exponent. For $v = (0,1)$, the norm of $\mathcal{A}(\omega, n)v$ grows linearly in $n$ and, hence, $v$ does not satisfy \eqref{eq:statement_lemma}. However, $w=(1,0)$ satisfies $\mathcal{A}(\omega, n)w=w$ and, hence, $w$ does fulfill~\eqref{eq:statement_lemma}. Note that $w$ is the 'terminal vector' of a Jordan block. Lemma \ref{lem:zero_LE_regularity} is closely related to the Jordan normal form of cocycles introduced in \cite{arnold1999jordan}. In the proof we used that $(\Omega, \mathcal{F})$ is a Polish space. However, we think that with slight adaptations to the proof, the lemma can be extended to more general probability spaces.
\end{remark}
A consequence of Lemma~\ref{lem:zero_LE_regularity} is the following result.
\begin{proposition}\label{prop:mane_sequence}
    Let $\mathcal{A}$ be a compact cocycle satisfying $\log^+ \|A(\cdot)\|\in L^1(\Omega, \mathbb R)$ and with a zero Lyapunov exponent.  Take $C>0$. For $\mathbb P$-a.e.~$\omega \in \Omega$ there are measurable  sequences $\{x(\omega, n)\}_{n\in \mathbb N}$ and $\{y(\omega, n)\}_{n\in \mathbb N}$ in $X$ (i.e. $\omega \mapsto x(\omega, n)$ and $\omega \mapsto y(\omega, n)$ are measurable maps for each $n\in \mathbb N$) and a measurable random variable $N(\omega)\in  \mathbb N$ such that the following conditions hold:
    \begin{enumerate}[(i)]
        \item $x(\omega, 0) = y(\omega, 0)$;
        \item $x(\omega, n+1) - A(\sigma^n \omega) x(\omega, n) = y(\omega, n+1)$, $\forall n\in \mathbb N$;
        \item $x(\omega, n) = y(\omega, n) = 0$ for $n>N(\omega)$;
        \item $\sup_{n\in \mathbb N} \| x(\omega, n)\| \geq C$ and $\sup_{n\in \mathbb N} \| y(\omega, n)\| \leq 1$.
        \end{enumerate}
\end{proposition}
\begin{proof}
    Let $\omega \mapsto v(\omega)$ be the measurable map with $\| v(\omega)\|=1$ given by  Lemma \ref{lem:zero_LE_regularity}. Fix $\omega \in \Omega$  such that the statement of the lemma, i.e.~equation \eqref{eq:statement_lemma}, holds.
    We define  sequences $\{x(\omega, n)\}_{n\in \mathbb N}$ and $\{y(\omega, n)\}_{n\in \mathbb N}$  by
    \begin{align}\label{eq:construction_x_y}
    \begin{split}
         y(\omega, n) &:= \beta(\omega, n) \|\mathcal{A}(\omega, n)v(\omega)\|^{-1} \mathcal{A}(\omega, n) v(\omega) ,\\
         x(\omega, n) &:= \sum_{m=0}^n \mathcal{A}(\sigma^m \omega, n-m)  y(\omega, m) =   \mathcal{A}(\omega, n) v(\omega) \sum_{m=0}^n \beta(\omega, m) \|\mathcal{A}(\omega, m)v(\omega)\|^{-1},
    \end{split}
    \end{align} 
    where $\beta(\omega, n) \in \mathbb R$.
    The sequences $\{ x(\omega, n)\}_{n\in \mathbb N}$ and $\{ y(\omega, n)\}_{n\in \mathbb N}$ are constructed such that $ x(\omega, 0) =  y(\omega, 0)= \beta(\omega, 0) v(\omega)$, 
    \[ x(\omega, n+1) - A(\sigma^n \omega)  x(\omega, n) = y(\omega, n+1),\]
    and $\| y(\omega, n)\| = |\beta(\omega, n)|$ for all $n\in \mathbb N$. Hence, we restrict ourselves to $|\beta(\omega, n)| \leq 1$. It remains to choose $\beta(\omega, n)$ measurably so that $\sup_{n\in \mathbb N}\|x(\omega,n)\|\geq C$ and such that both $x(\omega, n)$ and $y(\omega, n)$ are eventually zero.
    
    Note that $\|x(\omega, n)\|$ is given by $|\alpha(\omega, n)| \| \mathcal{A}(\omega, n) v(\omega)\|$, where
    \begin{equation*}
        \alpha(\omega, n) := \sum_{m=0}^n \beta(\omega, m)\|\mathcal{A}(\omega, m)v(\omega)\|^{-1}.
    \end{equation*}
    As long as $\beta(\omega, n) >0$, the sequence $\{\alpha(\omega, n)\}_{n\in \mathbb N}$ is increasing while for $\beta(\omega, n) <0$ the sequence decreases. Let $N^*(\omega), N(\omega) \in \mathbb N$ and define
    \begin{equation}\label{eq:beta_definition}
        \beta(\omega, n) := \begin{cases}
            +1, \quad  &\hspace{9mm}0\leq n \leq N^*(\omega),\\
            -1, \quad& N^*(\omega) < n < N(\omega), \\
            \beta(\omega, N(\omega)), \quad &\hspace{17mm}n=N(\omega), \\
            0,\quad &\hspace{17mm}n>N(\omega),
        \end{cases}
    \end{equation}
    where the value of $\beta(\omega, N(\omega))$ will be determined later.
    Our strategy is to measurably choose $N^*(\omega)$, $N(\omega)$ and $\beta(\omega, N(\omega))$ so that $\|x(\omega, N^*(\omega))\| \geq C$, and $\alpha(\omega, N(\omega) ) = 0$.
    
    The sequence $\{\alpha(\omega, n)\}_{n\in \mathbb N}$ grows monotonically for $n\leq N^*(\omega)$. By Lemma \ref{lem:zero_LE_regularity}, for $N^*(\omega)$ large enough, we find $\alpha(\omega, N^*(\omega))\geq 2C$. Additionally, Lemma \ref{lem:zero_LE_regularity} implies that there is a strictly increasing subsequence $\{n_k\}_{k\in \mathbb N}$ of $\mathbb N$ such that $\| x(\omega, n_k)\| \geq \frac 1 2 |\alpha(\omega, n_k)|$. Hence, for $N^*(\omega)$ large enough, we have $\|x(\omega, N^*(\omega))\| \geq C$. Let $N^*(\omega)\in \mathbb N$ be the smallest integer such that $\|  x(\omega, N^*(\omega))\| \geq C$. By construction, $\omega \mapsto N^*(\omega)$ is measurable.

    For $N^*(\omega) < n < N(\omega)$, the sequence $\{\alpha(\omega, n)\}_{n\in \mathbb N}$ decays monotonically. We have \[\alpha(\omega, n) = \alpha(\omega, N^*(\omega)) - \sum_{m=N^*(\omega) + 1}^n \|\mathcal{A}(\omega, m)v(\omega)\|^{-1}.\]
    By Lemma \ref{lem:zero_LE_regularity}, the sum on the right hand side grows indefinitely in $n$. Let $N(\omega)\in \mathbb N$ be the unique integer such that
    \[\sum_{m=N^*(\omega) + 1}^{N(\omega)-1} \|\mathcal{A}(\omega, m)v(\omega)\|^{-1} < \alpha(\omega, N^*(\omega)) \leq  \sum_{m=N^*(\omega) + 1}^{N(\omega)} \|\mathcal{A}(\omega, m)v(\omega)\|^{-1} .\]
    The map $\omega \mapsto N(\omega)$ is constructed in a measurable way. We now define
    \[\beta(\omega, N(\omega)):= -\left( \alpha(\omega, N^*(\omega)) - \sum_{m=N^*(\omega) + 1}^{N(\omega)-1} \|\mathcal{A}(\omega, m)v(\omega)\|^{-1} \right) \|\mathcal{A}(\omega, N(\omega) )v(\omega)\|.\]
    One can check that $|\beta(\omega, N(\omega)) | \leq 1$ and it is constructed so that $\alpha(\omega, N(\omega)) = 0$.

    We conclude that the sequences $\{x(\omega, n)\}_{n\in \mathbb N}$ and $\{y(\omega, n)\}_{n\in \mathbb N}$ as defined in \eqref{eq:construction_x_y} with coefficients $\beta(\omega, n)$ given by~\eqref{eq:beta_definition} have the desired properties.

\end{proof}

We are now in a position to give the proof of Theorem~\ref{converse}.
\begin{proof}[Proof of Theorem~\ref{converse}]
    We assume that the pair $(L^\infty(\Omega, X), L_C^\infty(\Omega, X))$ is admissible. For every $g\in L_C^\infty(\Omega, X)$ there is a unique $f \in L^\infty(\Omega, X)$ satisfying the admissibility condition \eqref{adm}. By Lemma \ref{lem:boundedness}, there is a constant $L>0$ such that any such pair satisfies $\|f\|_\infty \leq L\|g\|_{\infty, C}$. Assume for contradiction that there is no tempered exponential dichotomy. Then, by Proposition~\ref{criteria} there is a zero Lyapunov exponent for $\mathcal A$. We lead this assumption to a contradiction by constructing a pair $f \in L^\infty(\Omega, X)$ and $g\in L_C^\infty(\Omega, X)$ that satisfy the admissibility condition \eqref{adm}, but violating the bound of Lemma \ref{lem:boundedness}, i.e.~$\|f\|_\infty > L\|g\|_{\infty, C}$.
    
    Let $F\subset \Omega$ be a compact set with positive $\mathbb P$-measure on which $C$ is continuous. The existence of such a set is guaranteed by Lusin's theorem. Then, there is a constant $M>0$ such that $C(\omega) \le M$ for all $\omega \in F$. For $\omega \in F$, let 
    \begin{equation*}
        \tau_F(\omega) = \min \{m \geq 1 \mid \sigma^m\omega \in F\},
    \end{equation*}
    be the first return time to the set $F$ and let $\tau_F(\omega, n)$ be the time of the $n$-th return with the convention $\tau_F(\omega, 0) = 0$. By the Poincar\'e recurrence theorem, $\tau_F(\omega, n)$ is finite for $\mathbb P$-a.e.~$\omega \in F$. Let $\overline \sigma:F \to F$, $\omega \mapsto \sigma^{\tau_F(\omega)} \omega$ be the  return map of $\sigma$ to $F$, which is ergodic with respect to the conditional measure $\mathbb P_F := \mathbb P(F)^{-1} \mathbb P$ on $F$. Let $\overline{\mathcal{A}}$ be the induced cocycle over $\overline \sigma$ with generator $\overline{A}(\omega):=\mathcal{A}(\omega, \tau_F(\omega))$. Naturally, $\overline{\mathcal{A}}$ is a compact cocycle. From subadditivity of the logarithm and Kac's Theorem, e.g. \cite[Theorem 1.7]{sarig2009lecture}, we obtain
    \begin{equation*}
        \int_F \log^+\| \overline{A}(\omega)\| \, d\mathbb P_F \leq \int_F \sum_{n=0}^{\tau_F (\omega) -1}\log^+\| A(\sigma^n \omega)\|\, d\mathbb P_F = \frac{1}{\mathbb P(F)} \int_\Omega \log^+\| A(\omega) \|\,  d\mathbb P < \infty.
    \end{equation*}
    This shows that induced cocycle $\overline{\mathcal{A}}$ satisfies the assumptions of the multiplicative ergodic theorem. One can easily check that the Lyapunov exponents of $\overline{\mathcal{A}}$ are of the form $\mathbb P(F)^{-1}\lambda$, where $\lambda$ is a Lyapunov exponent of $\mathcal{A}$. In particular, since $\mathcal{A}$ has a zero Lyapunov exponent, $\overline{ \mathcal{A}}$ has a zero Lyapunov exponent as well. Applying Proposition \ref{prop:mane_sequence} to $\overline{\mathcal{A}}$, we obtain for $\mathbb P$-a.e.~$\omega$ in $F$ a measurable integer $N(\omega)$ and measurable sequences $\{\overline x(\omega, n)\}_{n\in \mathbb N}$ and $\{\overline y(\omega, n)\}_{n\in \mathbb N}$ that satisfy the following
    \begin{enumerate}[$(i)$]
        \item $\overline x(\omega, 0) = \overline y(\omega, 0)$;
        \item $\overline x(\omega, n+1) - \overline{A}(\overline \sigma^n \omega) \overline x(\omega, n) = \overline y(\omega, n+1)$, $\forall n\in \mathbb N$;
        \item $\overline x(\omega, n) = \overline y(\omega, n) = 0$ for $n$ greater than $N(\omega)$;
        \item $\sup_{n\in \mathbb N} \| \overline x(\omega, n)\| \geq (L+1)M $ and $\sup_{n\in \mathbb N} \| \overline y(\omega, n)\| \leq 1$.
    \end{enumerate}
    These sequences are constructed with respect to the induced cocycle $\overline{\mathcal{A}}$. We can extend these sequences to sequences $\{x(\omega, n)\}_{n\in \mathbb N}$, $\{y(\omega, n)\}_{n\in \mathbb N}$ with respect to $\mathcal{A}$. 
    For $\omega \in F$ and $n\in \mathbb N$ we find the unique $m\in \mathbb N$ such that $\tau_F(\omega, m)\le n<\tau_F(\omega, m+1)$ and set
    \[
    x(\omega, n):=\mathcal A(\sigma^{\tau_F(\omega, m)}\omega, n-\tau_F(\omega, m))\overline x(\omega, m)=\mathcal A(\overline{\sigma}^m \omega, n-\tau_F(\omega, m)\overline x(\omega, m).
    \]
    and 
    \[
    y(\omega, n)=\begin{cases}
  0 & \tau_F(\omega, m)<n; \\
  \overline y(\omega, m) & n=\tau_F(\omega, m).
    \end{cases}
    \]
    Note that $y(\omega, n)$ is non-zero only if $\sigma^n(\omega) \in F$. These sequences are 'interpolations' of $\{\overline x(\omega, n)\}_{n\in \mathbb N}$ and $\{\overline y(\omega, n)\}_{n\in \mathbb N}$ and inherit their properties, i.e.
    \begin{enumerate}[$(i)'$]
        \item $x(\omega, 0) = y(\omega, 0)$;
        \item $x(\omega, n+1) - A( \sigma^n \omega) x(\omega, n) = y(\omega, n+1)$, $\forall n\in \mathbb N$;
        \item $x(\omega, n) = y(\omega, n) = 0$ for $n$ greater than $\tau_F(\omega, N(\omega))$;
        \item $\sup_{n\in \mathbb N} \| x(\omega, n)\| \geq (L+1)M $ and $\sup_{n\in \mathbb N} \| y(\omega, n)\| \leq 1$.
    \end{enumerate}
    Except for statement $(ii)'$, these follow directly from the definition of $\{x(\omega, n)\}_{n\in \mathbb N}$ and $\{y(\omega, n)\}_{n\in \mathbb N}$.
    
    We verify $(ii)'$. Given $\omega \in F$ and $n\in \mathbb N$, let $m\in \mathbb N$ be as defined above. We first consider the case when $n+1<\tau_F(\omega, m+1)$. Then, 
    \[
    \begin{split}
x(\omega, n+1)-A(\sigma^n \omega)x(\omega, n) &=\mathcal A(\sigma^{\tau_F(\omega, m)}\omega, n+1-\tau_F(\omega, m))\overline x(\omega, m) \\
&\phantom{=}-A(\sigma^n \omega)\mathcal A(\sigma^{\tau_F(\omega, m)}\omega, n-\tau_F(\omega, m))\overline x(\omega, m) \\
&=\mathcal A(\sigma^{\tau_F(\omega, m)}\omega, n+1-\tau_F(\omega, m))\overline x(\omega, m) \\
&\phantom{=}-\mathcal A(\sigma^{\tau_F(\omega, m)}\omega, n+1-\tau_F(\omega, m))\overline x(\omega, m). \\
&=0=y(\omega, n+1).
    \end{split}
    \]
Let us now consider the case when $n+1=\tau_F(\omega, m+1)$.
Then,
\[
\begin{split}
x(\omega, n+1)-A(\sigma^n \omega)x(\omega, n) &=\overline x(\omega, m+1)-A(\sigma^n \omega)\mathcal A(\sigma^{\tau_F(\omega, m)}\omega, n-\tau_F(\omega, m))\overline x(\omega, m) \\
&=\overline x(\omega, m+1)-\mathcal A(\sigma^{\tau_F(\omega, m)}\omega, n+1-\tau_F(\omega, m))\overline x(\omega, m)\\
&=\overline x(\omega, m+1)-\mathcal A(\sigma^{\tau_F(\omega, m)}\omega, \tau_F(\omega, m+1)-\tau_F(\omega, m))\overline x(\omega, m)\\
&=\overline x(\omega, m+1)-\overline{A}(\overline \sigma^m\omega)\overline x(\omega, m) \\
&=\overline {y}(\omega, m+1)=y(\omega, n+1).
\end{split}
\]
    This shows that the sequences $\{x(\omega, n)\}_{n\in \mathbb N}$ and $\{y(\omega, n)\}_{n\in \mathbb N}$ satisfy $(ii)'$.

    Let $N \in \mathbb N$ be such that the set
    \begin{equation*}
        F_N := \{\omega \in F \mid \tau_F(\omega, N(\omega)) \leq N\},
    \end{equation*}
    has positive measure. By the Rokhlin Lemma, e.g.~\cite[Theorem 1.11]{jones1974transformations}, there is a measurable set $B\subset F_N$ with $\mathbb P(B)>0$ such that $\sigma^n(B) \neq \sigma^m(B)$ for all $n \neq m$ between $0$ and $N+1$.
    
    Define two function $f$, $g$ from $\Omega$ to $X$ by
    \begin{align*}
        f(\omega) &:= \sum_{n=0}^{N+1} \mathbf{1}_{\sigma^n B}(\omega) x(\sigma^{-n} \omega, n),\\
        g(\omega) &:= \sum_{n=0}^{N+1} \mathbf{1}_{\sigma^n B}(\omega) y(\sigma^{-n} \omega, n).
    \end{align*}
    We write the sums until $N+1$, but the last summand is zero since $x(\omega, N+1)=y(\omega, N+1) = 0$. Hence, $f$ and $g$ are supported on $\cup_{n=0}^N \sigma^n B$. 
    We verify that $f$ and $g$ form an admissible pair, i.e.~that~\eqref{adm} holds. We distinguish three cases. For $\omega \in B$, we have that $f(\omega)=x(\omega, 0)$, $f(\sigma^{-1}\omega)=0$ and $g(\omega)=y(\omega, 0)=x(\omega, 0)$. Consequently, 
    \[
    f(\omega)-A(\sigma^{-1}\omega)f(\sigma^{-1}\omega)=x(\omega, 0)=g(\omega).
    \]
    We now consider the case when $\omega \in \sigma^k B$ for some $1\le k\le N+1$. Observe that $\sigma^{-1}\omega \in \sigma^{k-1}B$. Hence, 
 $f(\omega)=x(\sigma^{-k}\omega, k)$, $f(\sigma^{-1}\omega)=x(\sigma^{-k}\omega, k-1)$ and 
 $g(\omega)=y(\sigma^{-k}\omega, k)$. By $(ii)'$ we have that 
 \[
 f(\omega)- A(\sigma^{-1}\omega)f(\sigma^{-1}\omega)=x(\sigma^{-k}\omega, k)- A(\sigma^{-1}\omega)x(\sigma^{-k}\omega, k-1)=y(\sigma^{-k}\omega, k)=g(\omega).
 \]
    Finally, in the case when $\omega \notin \sigma^k B$ for $0\le k\le N+1$, then $\sigma^{-1}\omega \notin \sigma^k B$ for $0\le k\le N$. Consequently, 
    $f(\omega)=f(\sigma^{-1}\omega)=g(\omega)=0$, and thus
    $
     f(\omega)- A(\sigma^{-1}\omega)f(\sigma^{-1}\omega)=g(\omega).
    $
    This completes the proof that $f$ and $g$ satisfy the admissibility condition \eqref{adm}.

    For each $\omega \in B$, there is an $0\leq n \leq N$ for which $\|x(\omega, n)\| \geq (L+1)M$. Hence, there is an $0\leq n \leq N$ such that $\|x(\omega, n)\| \geq (L+1)M$ for all $\omega$ in a subset of $B$ of positive measure. We conclude that the essential supremum of $f$ on $\sigma^nB$ is at least $(L+1)M$, and thereby $\| f\|_\infty \geq (L+1)M$.

    Lastly, observe that $g(\omega)$ is only non-zero if $\omega \in \sigma^n B$ and $y(\sigma^{-n}\omega, n)\neq 0$. But $y(\sigma^{-n}\omega, n)$ is only non-zero if $\sigma^n \circ \sigma^{-n} \omega = \omega$ is in $F$. Hence, $g$ is supported on $F$. Since $\|y(\omega, n)\|\leq 1$, we have $\|g(\omega)\| \leq 1$ for all $\omega \in \Omega$ and conclude $\|g\|_{\infty, C} \leq M$. 

    This shows $\|f\|_\infty > LM \geq L\|g\|_{\infty, C}$ which contradicts Lemma \ref{lem:boundedness}. This shows that $\mathcal{A}$ must admit a tempered exponential dichotomy.
\end{proof}

The proof of Theorem \ref{converse2} is similar to the proof of Theorem \ref{converse} and is therefore not given in full detail.
\begin{proof}[Proof of Theorem \ref{converse2}]
    Let $L>0$ be the constant from Lemma \ref{lem:boundedness} such that any admissible pair of functions $f\in L_C^\infty(\Omega, X)$ and $g \in L^\infty(\Omega, X)$ satisfies $\|f\|_{\infty,C} \leq L\|g\|_{\infty}$. Assume for contradiction that there is no tempered exponential dichotomy. Then, by Proposition~\ref{criteria} there is a zero Lyapunov exponent for $\mathcal A$. We lead this assumption to a contradiction by constructing an admissible pair of functions that violate the bound of Lemma \ref{lem:boundedness}, i.e.~$\|f\|_{\infty, C} > L\|g\|_\infty$.
    
    Let $F\subset \Omega$ be a compact set with positive $\mathbb P$-measure on which $C$ is continuous. There is a constant $M>0$ such that $C(\omega) \ge M$ for all $\omega \in F$. Let $\overline \sigma:F \to F$, $\omega \mapsto \sigma^{\tau_F(\omega)} \omega$ be the return map of $\sigma$ to $F$ and let $\overline{\mathcal{A}}$ be the induced cocycle over $\overline \sigma$ with generator $\overline{A}$. The induced cocycle $\overline{ \mathcal{A}}$ over $F$ still has a zero Lyapunov exponent. Applying Proposition \ref{prop:mane_sequence} to $\overline{\mathcal{A}}$, we obtain for $\mathbb P$-a.e.~$\omega$ in $F$ a measurable integer $N(\omega)$ and measurable sequences $\{\overline x(\omega, n)\}_{n\in \mathbb N}$ and $\{\overline y(\omega, n)\}_{n\in \mathbb N}$ that satisfy
    \begin{enumerate}[$(i)$]
        \item $\overline x(\omega, 0) = \overline y(\omega, 0)$;
        \item $\overline x(\omega, n+1) - \overline{A}(\overline \sigma^n \omega) \overline x(\omega, n) = \overline y(\omega, n+1)$, $\forall n\in \mathbb N$;
        \item $\overline x(\omega, n) = \overline y(\omega, n) = 0$ for $n$ greater than $N(\omega)$;
        \item $\sup_{n\in \mathbb N} \| \overline x(\omega, n)\| \geq (L+1)M^{-1} $ and $\sup_{n\in \mathbb N} \| \overline y(\omega, n)\| \leq 1$.
    \end{enumerate}
    These sequences are constructed with respect to the induced cocycle $\overline{\mathcal{A}}$. Like in the proof of Theorem \ref{converse}, we extend these sequences to sequences $\{x(\omega, n)\}_{n\in \mathbb N}$, $\{y(\omega, n)\}_{n\in \mathbb N}$ with respect to $\mathcal{A}$, where $\omega \in F$.
    These sequences satisfy
    \begin{enumerate}[$(i)'$]
        \item $x(\omega, 0) = y(\omega, 0)$;
        \item $x(\omega, n+1) - A( \sigma^n \omega, 1) x(\omega, n) = y(\omega, n+1)$, $\forall n\in \mathbb N$;
        \item $x(\omega, n) = y(\omega, n) = 0$ for $n$ greater than $\tau_F(\omega, N(\omega))$;
        \item $\sup_{n\in \mathbb N} \| x(\omega, n)\| \geq (L+1)M^{-1} $ and $\sup_{n\in \mathbb N} \| y(\omega, n)\| \leq 1$.
    \end{enumerate}
    Additionally, for $\mathbb P$-a.e.~$\omega \in F$, there is an $n\in \mathbb N$ such that $\| x(\omega, n)\| \geq (L+1)M^{-1}$ and $\sigma^n \omega \in F$.

    Let $N \in \mathbb N$ be such that the set
    \begin{equation*}
        F_N := \{\omega \in F \mid \tau_F(\omega, N(\omega)) \leq N\},
    \end{equation*}
    has positive measure. By the Rokhlin Lemma, 
    there is a measurable set $B\subset F_N$ with $\mathbb P(B)>0$ such that $\sigma^n(B) \neq \sigma^m(B)$ for all $n \neq m$ between $0$ and $N+1$.
    
    Define two function $f$, $g$ from $\Omega$ to $X$ by
    \begin{align*}
        f(\omega) &:= \sum_{n=0}^{N+1} \mathbf{1}_{\sigma^n B}(\omega) x(\sigma^{-n} \omega, n),\\
        g(\omega) &:= \sum_{n=0}^{N+1} \mathbf{1}_{\sigma^n B}(\omega) y(\sigma^{-n} \omega, n).
    \end{align*}
    The proof that $f$ and $g$ satisfy the admissibility condition \eqref{adm} is the same as in the proof of Theorem \ref{converse2}.

    For each $\omega \in B$, there is an $0\leq n \leq N$ such that $\|x(\omega, n)\| \geq (L+1)M^{-1}$ and such that $\sigma^n \omega \in F$. Hence, there is an $0\leq n \leq N$ such that $\|x(\omega, n)\| \geq (L+1)M^{-1}$ and $\sigma^n \omega \in F$ for all $\omega$ in a subset of $B$ of positive measure. We conclude $\sigma^n B\cap F$ is a set of positive measure and that the essential supremum of $f$ on this set is at least $(L+1)M^{-1}$. We conclude $\| f\|_{\infty,C} \geq L+1$.

    Lastly, observe that by construction $\|y(\omega, n)\|\leq 1$, and, thereby, $\|g\|_{\infty} \leq 1$. This shows $\|f\|_{\infty,C} > L \geq L\|g\|_{\infty}$ which contradicts Lemma \ref{lem:boundedness}. This shows that $\mathcal{A}$ must admit a tempered exponential dichotomy.
\end{proof} 

\begin{remark}\label{remark}
    The assumption that the linear cocycle $\mathcal{A}$ is compact is only used for the applicability of the multiplicative ergodic theorem (Theorem \ref{MET}), which allows us to characterize the existence of a tempered exponential dichotomy in terms of the existence of a zero Lyapunov exponent (see Proposition \ref{criteria}). The multiplicative ergodic theorem from \cite{gonzalez2014semi} holds not only for compact cocycles, but more generally for quasi-compact cocycles. Additionally, Proposition \ref{criteria} holds for quasi-compact cocycles as well, as long as the index of compactness is less then 0, cf.~\cite[Proposition 3.2]{backes2019periodic}. Hence, the proofs of this section work without modification for quasi-compact cocycles with an index of compactness of less than $0$. 
\end{remark}

\begin{remark}
In~\cite{dragivcevic2023one} the authors have obtained a Chow-Leiva type characterization of tempered exponential dichotomies in which the input spaces consists of bounded sequences in $X$ and the output space of weighted-bounded sequences in $X$, which is close in spirit to Propositions~\ref{thm:L_LC} and Theorem~\ref{converse}. We do not see how to use~\cite{dragivcevic2023one} to obtain the results of the present paper since the strategy outlined in~\cite{latushkin1999spectral} relies on the fact that in their setting the input and output space coincide and the Mather operator is a well-defined bounded linear operator. 

For other results dealing with characterizations of tempered dichotomies using Chow-Leiva approach, we refer to~\cite{abu2017spectral,barreira2016tempered,  dragivcevic2011characterization,dragivcevic2024measurable}.

\end{remark}

\section{Robustness of tempered exponential dichotomies}\label{sec:robustness}
We apply our results to establish the persistence of the notion of a tempered exponential dichotomy under small (linear) perturbations. We stress that this application is rather standard and that the robustness property of tempered dichotomies was established for cocycles which do not necessarily satisfy assumptions of the Oseledets theorem in~\cite{zhou2013roughness} (although with a much longer proof).

\begin{theorem}
Let $\mathcal A\colon \Omega \times \mathbb N\to \mathcal L(X)$ be a linear cocycle that admits a tempered exponential dichotomy. Then, there exists a tempered random variable $c\colon \Omega \to (0, \infty)$ with the property that each compact cocycle $\mathcal B\colon \Omega \times \mathbb N\to \mathcal L(X)$ whose generator $B$ satisfies
\begin{equation}\label{rob}
\|A(\omega)-B(\omega)\| \le c(\omega) \quad \text{for $\mathbb P$-a.e. $\omega \in \Omega$}
\end{equation}
also admits a tempered exponential dichotomy.
\end{theorem}

\begin{proof}
Let $K\colon \Omega \to (0, \infty)$ and $\lambda>0$ be as in Definition~\ref{TED} (associated to the tempered exponential dichotomy of $\mathcal A$). We choose $d>0$ such that
\begin{equation}\label{d}
d\cdot \frac{1+e^{-\lambda}}{1-e^{-\lambda}}<1,
\end{equation}
and we set \begin{equation}\label{c} c(\omega):=\frac{d}{K(\sigma \omega)}, \quad  \omega \in \Omega.
\end{equation}Since $K$ is tempered we have that $c$
is also tempered.

Fix an arbitrary $g\in L_K^\infty(\Omega, X)$. We define $\mathcal T\colon L^\infty(\Omega, X)\to L^\infty(\Omega, X)$ by 
\[
\begin{split}
(\mathcal T f)(\omega) &=\sum_{n=0}^\infty \mathcal A(\sigma^{-n}\omega, n)\Pi^s(\sigma^{-n}\omega)\left ((B(\sigma^{-(n+1)}\omega)-A(\sigma^{-(n+1)}\omega))f(\sigma^{-(n+1)}\omega)+g(\sigma^{-n}\omega) \right ) \\
&\phantom{=}-\sum_{n=1}^\infty \mathcal A(\sigma^n\omega, -n)\Pi^u(\sigma^n \omega) \left ((B(\sigma^{n-1}\omega)-A(\sigma^{n-1}\omega))f(\sigma^{n-1}\omega)+g(\sigma^n \omega)\right ),
\end{split}
\]
for $\omega \in \Omega$ and $f\in L^\infty(\Omega, X)$.  Take $f\in L^\infty (\Omega, X)$. By~\eqref{dic1}, \eqref{dic2}, \eqref{rob} and~\eqref{c} we have that 
\[
\begin{split}
\|(\mathcal T f)(\omega)\| &\le \sum_{n=0}^\infty K(\sigma^{-n}\omega)e^{-\lambda n}(c(\sigma^{-(n+1)}\omega)\|f(\sigma^{-(n+1)}\omega)\|+\|g(\sigma^{-n}\omega)\|) \\
&\phantom{\le}+\sum_{n=1}^\infty K(\sigma^n \omega)e^{-\lambda n}(c(\sigma^{n-1}\omega)\|f(\sigma^{n-1}\omega)\|+g(\sigma^n \omega)\|) \\
&\le \left (d\|f\|_\infty+\|g\|_{\infty, K}\right)\sum_{n=0}^\infty e^{-\lambda n}+\left (d\|f\|_\infty+\|g\|_{\infty, K}\right)\sum_{n=1}^\infty e^{-\lambda n} \\
&=\frac{1+e^{-\lambda}}{1-e^{-\lambda}}\left (d\|f\|_\infty+\|g\|_{\infty, K}\right),
\end{split}
\]
for $\mathbb P$-a.e. $\omega \in \Omega$. Hence, 
\[
\|\mathcal Tf\|_\infty \le \frac{1+e^{-\lambda}}{1-e^{-\lambda}}\left (d\|f\|_\infty+\|g\|_{\infty, K}\right),
\]
which implies that $\mathcal Tf \in L^\infty(\Omega, X)$. Thus, $\mathcal T$ is well-defined. Similarly,  we find that 
\begin{equation}\label{contraction}
\|\mathcal Tf_1-\mathcal T f_2\|_\infty \le d\cdot \frac{1+e^{-\lambda}}{1-e^{-\lambda}}\|f_1-f_2\|_\infty, \quad \text{for $f_1, f_2\in L^\infty(\Omega, X)$.}
\end{equation}
Due to our choice of $d$ (see~\eqref{d}), we conclude that $\mathcal T$ is a contraction of $L^\infty(\Omega, X)$, and therefore it has a unique fixed point $f\in L^\infty(\Omega, X)$.  Hence,
\[
\begin{split}
f(\omega) &=\sum_{n=0}^\infty \mathcal A(\sigma^{-n}\omega, n)\Pi^s(\sigma^{-n}\omega)\left ((B(\sigma^{-(n+1)}\omega)-A(\sigma^{-(n+1)}\omega))f(\sigma^{-(n+1)}\omega)+g(\sigma^{-n}\omega) \right ) \\
&\phantom{=}-\sum_{n=1}^\infty \mathcal A(\sigma^n\omega, -n)\Pi^u(\sigma^n \omega) \left ((B(\sigma^{n-1}\omega)-A(\sigma^{n-1}\omega))f(\sigma^{n-1}\omega)+g(\sigma^n \omega)\right ),
\end{split}
\]
for $\mathbb P$-a.e. $\omega \in \Omega$. A computation analogous to that in~\eqref{computation} shows that 
\[
f(\omega)-A(\sigma^{-1}\omega)f(\sigma^{-1}\omega)=(B(\sigma^{-1}\omega)-A(\sigma^{-1}\omega))f(\sigma^{-1}\omega)+g(\omega),
\]
for $\mathbb P$-a.e. $\omega \in \Omega$. In other words,
\[
f(\omega)-B(\sigma^{-1}\omega)f(\sigma^{-1}\omega)=g(\omega), \quad \text{for $\mathbb P$-a.e. $\omega \in \Omega$.}
\]
Therefore, $f$ and $g$ are an admissible pair with respect to the cocycle $\mathcal B$. The uniqueness of $f$ follows from the uniqueness of a fixed point for $\mathcal T$. We conclude that the pair $(L^\infty(\Omega, X), L_K^\infty(\Omega, X))$ is admissible for $\mathcal B$, and consequently Theorem~\ref{converse} implies that $\mathcal B$ admits a tempered exponential dichotomy. The proof of the theorem is completed.
\end{proof}

\section*{Acknowledgements}
D.~Dragi\v cevi\' c was supported in part by University of Rijeka under the project uniri-iskusni-prirod-23-98
3046.

\bibliographystyle{alpha}
\bibliography{bibliography.bib}

\end{document}